%% file: Luchko-Mainardi_newJVA-E-print.tex
\documentclass[11pt,a4paper]{amsart}
\usepackage{graphicx}  
\numberwithin{equation}{section}
\allowdisplaybreaks
\begin{document}
\title[The time-fractional diffusion-wave equation]{Cauchy and signaling problems  for \\ the time-fractional diffusion-wave equation}
%: a survey of some known and novel  results
%%% first author
 \author[Yuri Luchko,  Francesco Mainardi]  {Yuri Luchko $^1$, Francesco Mainardi $^2$}
%\author{Yuri Luchko$^1$}
\address{${}^1$  Professor of Mathematics \\
	Department of Mathematics, Physics, and Chemistry, 
   Beuth Technical University of Applied Sciences, 
13353 Berlin, Germany.}
\email{luchko@beuth-hochschule.de}

\address{${}^2$   Professor of Mathematical Physics \\
Department of Physics and Astronomy,
Bologna University, and INFN,
40126 Bologna, Italy.}
\email{francesco.mainardi@unibo.it; \ francesco.mainardi@bo.infn.it }

\date{September 2016 \\
This E-Print is a reproduction with a different layout
 of the paper published  in\\
{\bf ASME Journal of Vibration and Acoustics,
Vol. 136 (2014), pp. 051008/1-7.} \\ DOI: 10.1115/1.4026892}

%\begin{document}

%\maketitle    

%%%%%%%%%%%%%%%%%%%%%%%%%%%%%%%%%%%%%%%%%%%%%%%%%%%%%%%%%%%%%%%%%%%%%%
\begin{abstract}
{In this paper, some known and novel properties of the Cauchy and signaling problems for the one-dimensional time-fractional diffusion-wave equation with the Caputo fractional derivative of order $\beta,\ 1 \le \beta \le 2$ are investigated. In particular, their response to a localized disturbance of the initial data is studied.  It is known that whereas the diffusion equation describes a process where the disturbance spreads infinitely fast, the propagation velocity of the disturbance is a constant for the wave equation.
We show that the time-fractional diffusion-wave equation interpolates between these two different responses in the sense that the propagation velocities of the maximum points, centers of gravity, and medians of  the fundamental solutions to both the Cauchy and the signaling problems are all finite. On the other hand, the disturbance spreads infinitely fast and the time-fractional diffusion-wave equation is non-relativistic like the classical diffusion equation.
In this paper, the maximum locations, the centers of gravity, and the medians of the fundamental solution to the Cauchy and signaling problems and their propagation velocities   are described analytically and calculated numerically. The obtained results for the Cauchy and the signaling problems are interpreted and compared to each other.}
\end{abstract}
\maketitle

%%%%%%%%%% MAINARDI DEFINITIONS
%%%%%%%%%%%%%%%%%%%%%%%%%%%%%%%%%%%%%%%%%%%%%%%%%%%%%%%%%%%%
\def\pni{\par \noindent}
\def\vsh{\smallskip}
\def\vs{\medskip}
\def\vvs{\bigskip}
\def\vvvs{\bigskip\medskip} %% {\vskip 1.5truecm}
\def\vsp{\par}
\def\vsn{\vsh\pni}
\def\cen{\centerline}
\def\ra{\item{a)\ }} \def\rb{\item{b)\ }}   \def\rc{\item{c)\ }}
\def\eg{{\it e.g.}\ } \def\ie{{\it i.e.}\ }
%% REFERENCES
\def\rstar{\item{$\null ^*$)\ }}
   \def\rl{\item{-- \ }}
    \def\rp{\item{}}
    \def\rf#1{\item{{#1}.}}  %% Please note N. !!!
%%%% MATHEMATICS
\def\e{{\rm e}}
\def\v{{\rm v}}
\def\V{{\rm V}}
\def\exp{{\rm exp}}
\def\ds{\displaystyle}
\def\dis{\displaystyle}
\def\q{\quad}	 \def\qq{\qquad}
\def\lan{\langle}\def\ran{\rangle}
\def\l{\left} \def\r{\right}
\def\lra{\Longleftrightarrow}
\def\d{\partial}
 \def\dr{\partial r}  \def\dt{\partial t}
\def\dx{\partial x}   \def\dy{\partial y}  \def\dz{\partial z}
\def\rec#1{{1\over{#1}}}
%%%%%%%%% for LAPLACE TRANSFORMS %%%%%%%%
\def\bar{\tilde}
\def\barr{\widetilde}
\def\epsilons{{\tilde \epsilon(s)}}
\def\sigmas{{\tilde \sigma (s)}}
\def\fs{{\tilde f(s)}}
\def\Js{{\tilde J(s)}}
\def\Gs{{\tilde G(s)}}
\def\Fs{{\tilde F(s)}}
 \def\Ls{{\tilde L(s)}}
\def\L{{\mathcal L}} %%% Laplace Transform !!!!
\def\F{{\mathcal F}} %%% Fourier Transform !!!!
%%%%%%%%%%%% SETS of NATURAL, REAL, COMPLEX NUMBERS : \NN, \RR, \CC
\def\NN{{\rm I\hskip-2pt N}}
\def\RR{\vbox {\hbox to 8.9pt {I\hskip-2.1pt R\hfil}}}
\def\CC{{\rm C\hskip-4.8pt \vrule height 6pt width 12000sp\hskip 5pt}}
\def\II{{\rm I\hskip-2pt I}}
%%%%%%%%%%%%
 \def\D{{\mathcal D}}
%%%%%%%%% GREEN FUNCTIONS %%%%%%%%
\def\Gc{{\mathcal {G}}_c}	\def\Gcs{\barr{\Gc}} %% CAUCHY PROBLEM
\def\Gs{{\mathcal {G}}_s}	\def\Gss{\barr{\Gs}} %% SIGNALLING PROBLEM
\def\args{x\, s^{1/2}}
\def\argsa{x/\, s^{\nu }}
\def\arg{ x^2/ (4 \, t)}
%%%%%%%%%%%%%%%%%%%%
%% CISM \line{4. THE FRACTIONAL DIFFUSION-WAVE EQUATION \hfill}
\section {Introduction}
By the {\it fractional diffusion-wave} equation we mean a  linear
integro partial differential equation
obtained from the
classical diffusion or wave equation by replacing the first- or
second-order time derivative by a fractional derivative (in the Caputo
sense) of order
$\beta	  $ with $0 <\beta    \le 2$. In our notations it reads
\begin{equation}
{\d^{\beta  }  u\over \dt^{\beta   }} =
 \D \,{\d^2 u \over \dx^2} \,,\qq u=u(x,t)\,, \q
0<\beta  \le 2 \,, \q \D>0\,,
\label{eqno(1)}
\end{equation}
where
$\D$ denotes a positive constant with the dimension $L^2\,T^{-\beta  }\,,$
$x$ and $t$ are the space and time variables,
and  $u=u(x,t)$ is the field variable, which is assumed
to be a {\it causal} function of time, \ie
vanishing for $t<0\,. $
 \vsp
 Recalling  the definition of the Caputo fractional derivative, see  e.g. 
Gorenflo and Mainardi \cite{Gorenflo-Mainardi_CISM1997},
Podlubny \cite{Podlubny_BOOK1999}, 
and setting for convenience, but  without loss of generality,   $\D \equiv 1$,
we get in explicit form the following integro-differential equations 
%% In Eq. (1) we thus need to distinguish  two cases
%%  $i)\q 0<\beta	\le 1\,, \q {\rm and}\q ii)\q 1<\beta  \le 2\,, $
%% for which the equation assumes the explicit forms as follows :
\begin{equation}
%%  \Phi_{1-\beta } (t) \,*\, {\d u\over \dt} =
    \rec{\Gamma(1-\beta )}\,
 \int_0^t (t-\tau )^{-\beta }\, \l( {\d u\over \d\tau}\r) \, d\tau
   =   {\d^2 u \over \dx^2} \,,\qq 0<\beta	\le 1 \,; 
   \label{eqno(2)}
   \end{equation}
\begin{equation}
%% \Phi_{2-\beta } (t) \,*\, {\d^2 u\over \dt^2} =
    \rec{\Gamma(2-\beta )}\,
 \int_0^t (t-\tau )^{1-\beta } \, \left({\d^2 u\over\d\tau^2}\right)\, d\tau
   =   {\d^2 u \over \dx^2} \,,\q 1<\beta  \le 2 \,. 
   \label{eqno(3)}
   \end{equation}
The equations (\ref{eqno(2)}) and (\ref{eqno(3)}) can be properly referred to as  the
{\it time-fractional diffusion} and the {\it time-fractional wave}
equation, respectively.
%%%%%%%%
\vsp
%%%%%%%
A {\it fractional diffusion} equation akin to (\ref{eqno(2)})
has been formerly introduced in 1986 by Nigmatullin \cite{Nigmatullin_1986}
to describe  diffusion in special types of porous media,
which exhibit a  fractal geometry.
In 1995 Mainardi \cite{Mainardi_1995a}  has shown  that
the {\it fractional wave} equation (\ref{eqno(3)}) governs the propagation of mechanical diffusive waves in viscoelastic media which exhibit a simple power-law
creep. This problem of dynamic viscoelasticity
 was formerly treated by Pipkin \cite{Pipkin_BOOK1986}
  and Kreiss and  Pipkin \cite{Kreis-Pipkin_1986}
who however were unaware of the interpretation by fractional calculus.
In our opinion, the above references  provide  some interesting 
and pioneering examples of the relevance of  (\ref{eqno(1)}) in physics.
Of course, any time some  hereditary mechanisms of power-law
type are present in  diffusion or wave phenomena, an
appearance of time-fractional derivatives in the evolution
equations  is expected.
\vsp
In a series of papers \cite{Mainardi_1994,Mainardi_1995a,Mainardi_1995b,%
Mainardi_1996a,Mainardi_1996b,Mainardi_CISM1997}, 
  Mainardi  has pursued his analysis
on the time-fractional diffusion-wave equation (\ref{eqno(1)})
%% in order to provide more insight
%% into fractional diffusion and fractional propagation.
based on Laplace transforms and special functions of Wright type.
%% has provided analytical and numerical results on
%% the fundamental solutions which turn out to be expressed
%% in terms of a special function of Wright type in the similarity
%% variable.
Other mathematical aspects of  integro-differential equations
akin to  (\ref{eqno(1)})  based on the use of the integral transforms
and special functions have been also treated in some relevant papers including 
those
by Wyss \cite{Wyss_JMP1986}, 
Schneider and  Wyss \cite{Schneider-Wyss_JMP1989}, 
%% (Mellin transforms and Fox $H$ functions)
Fujita \cite{Fujita_1990a},
 %(Fourier transforms and Mittag-Leffler functions).
 Pr\"uss \cite{Prusse_BOOK93},
Mainardi, Luchko and Pagnini \cite{Mainardi-Luchko-Pagnini_FCAA2001},
Mainardi, Pagnini and Saxena \cite{Mainardi-Pagnini-Saxena_2005},
and more recently by 
Luchko \cite{Luchko_2010,Luchko_2012,Luchko_2013}.
%More formal  approaches based on semigroup theory in Banach spaces
%have been given by Kochubei \cite{Kochubei_1989,Kochubei_1990}
%  and by El-Sayed \cite{El-Sayed_1996}.
%\vsp
Furthermore, mathematical aspects related to similarity properties and stable probability densities have been treated  
 by  Fujita \cite{Fujita_1990b},  Engler \cite{Engler_1997},
Mainardi and Tomirotti \cite{Mainardi-Tomirotti_GEO1997},
Luchko and Gorenflo \cite{Luchko_1998},
Gorenflo, Luchko and Mainardi \cite{GoLuMa_1999,GoLuMa_2000},
and more recently by 
Luchko, Mainardi and Povstenko \cite{LMP_2013}, 
and by Luchko and Mainardi \cite{Luchko-Mainardi_2013}. 
%% in a more abstract form
We also outline the papers \cite{Luchko_2009,Luchko_2011}
by Luchko on the application of the maximum principle 
to the time-fractional diffusion equations. 
Of course, the above  list of references is  not exhaustive and mainly regards those
 that  have attracted our attention.
\vsp
In this paper,   we consider Eq. (\ref{eqno(1)})
restricting out analysis to the case $1\le \beta\le  2$ 
that we refer to as {\it the time-fractional diffusion-wave equation}. 
Our main purpose is  to point out 
some relevant properties of the related intermediate process that governs  transition from pure diffusion ($\beta=1$) to pure wave propagation ($\beta=2$).
%% referring to the appendix
%% with the addition of some new considerations on the relevant connection
%% with the stable   probability distributions.
\vsp
In the second section,  we  define
the two basic boundary-value problems,
referred to as
the {\it Cauchy problem}  and the {\it Signaling problem},
recalling for them 
the respective fundamental solutions (the {\it Green
functions}).
We outline a reciprocity relation between the Green functions
themselves in the space-time domain.  In view of this relation
the Green functions can be expressed
in terms  of two interrelated {\it auxiliary functions}  in
the similarity variable $r = |x|/t^{\beta  /2}\,.$
In some plots, the evolution of the fundamental solutions
of both the Cauchy and Signaling problems for some values
of the order of the time derivative is shown. These solutions exhibit a  pulse-like pattern moving along the $x$ axis that depends of the order of the time-fractional  derivative.   
This allows us to better recognize
 the   
processes  intermediate between diffusion and wave propagation.
\vsp
In the third section, 
we  analyze the location and the evolution  of the maximum  
of the pulse like patterns  and its dependence  of the order of the fractional derivative.
In particular, we  present both an analytical treatment of the maximum locations, maximum values, and the propagation velocities of the maximum points of the Green functions and their plots. 
\vsp
Then in the fourth section we consider the location of the center of gravity and of the median for these pulse-like patterns in order to compare their evolution with respect to  that of the corresponding maximum.
\vsp 
 Finally, in the last section some conclusive remarks are given.

%%%%%%%%%%\vfill\eject%%%%%%%%%
\section{Cauchy and signaling problems}
As it is well known, the two basic boundary-value problems
for the  evolution equations of diffusion
and wave  type	are the {\it Cauchy} and the {\it Signaling problems}.
%In the	{\it Cauchy problem}, which concerns  the space-time domain
%$-\infty <x< + \infty\,, $ $\, t \ge 0\,, $
%the data are assigned at $t=0^+$ on
%the  whole space axis (initial data).
%In the {\it Signaling problem}, which concerns  the space-time domain
%$x\ge 0\,, $ $\, t \ge 0\,, $
%the data are assigned both at $t=0^+$ on the semi-infinite
%space axis $ x >0 $ (initial data) and at $x=0^+$ on the semi-infinite
%time axis  $ t>0$ (boundary data); here, as mostly usual,
%the initial data are assumed to be vanishing.
%%%%%%% VEDERE KEVORKIAN for the Signalling Problem !!!!
%% \vsp
Extending
the classical analysis to our fractional equation (\ref{eqno(1)}), and
denoting by $f(x)$ and $h(t)$ two given, sufficiently well-behaved
functions,    the basic problems are thus formulated as following:
\vsh\pni
\begin{equation}
\label{eqno(4a)}
% \left\{
\begin{array}{ll}
{\hbox { Cauchy problem:}} 
&  \; u(x,0^+)=f(x) , \; -\infty <x < +\infty,
\\  
& u(\mp \infty,t) = 0,\;  t>0;
\end{array}
\end{equation} 
%%%%%%%%%%    
  \begin{equation}
  \label{eqno(4b)}
  \begin{array}{ll}
%%%%%%%%%%
{\hbox   {Signaling problem:}} 
& \; u(x, 0^+) =0 , \; 0 <  x < +\infty,
\\
&    u(0^+,t ) =h(t), \; u(+\infty,t) =0 , \;   t >0 . 
\end{array}
%\right.
%\label{eqno(4a-4b}
\end{equation}
%%%%%% SI POTREBBERO mettere da parte le RADIATION CONDITIONS:
%%% VANISHING LIMITS for x \to \pm \infty or for x \tp +\infty
%% at any time %%%%%%%%%%%%%%%%%%%%%%%%%%%%%%%%%%%%%%%%
%% \vsp
For $1<\beta    \le 2\,, $  
the initial value of the first-order  time derivative 
of the field variable,
\ie ${\ds \frac{\d }{\d t} u(x,0^+)= g(x)}$,   is required  in the above problems, since in this case
the second time derivative   appears  in the integro-differential equation
 (\ref{eqno(3)}) and, consequently, two linearly independent solutions are to be
determined.   In what follows, we mainly limit ourselves to the case
$g(x) \equiv 0$.
\vsp
In view of our analysis, we find it convenient to put
%% introduce the {\it similarity variable}
%% $$ z := {|x|\over \sqrt{D} \, t^\nu } \,, \q
\begin{equation}
 \nu  ={\beta  \over 2}\,,
   \q 0<\nu  <1\,. 
\end{equation}
%%%%%% \vfill\eject
%%%%%%%%%
\vsp
For the {\it Cauchy} and {\it Signaling}  problems, we
introduce
the so-called Green functions $\Gc (x,t;\nu  )$ and $\Gs(x,t;\nu  )$,
which represent the fundamental solutions that are 
obtained when $f(x) = \delta (x)$ and $h(t) = \delta (t) $,
where $\delta$ denotes the Dirac $\delta$-function.
It should be noted that  
the Green function for the Cauchy problem turns out to be
an even function of $x$, so $\Gc(x ,t ;\nu ) = \Gc(|x|,t ;\nu )$.

For $0<\nu  \le 1$, the two Green functions 
%% First we note the  for the  Green functions
are connected by the following {\it reciprocity relation} (see the already cited papers by Mainardi): 
   \begin{equation}
 2\nu  \, x\,  \Gc(x,t;\nu  )  = t\, \Gs(x,t;\nu  )
    = F_\nu(r) = \nu  r\, M_\nu(r),\
  \label{eqno(15)} 
  \end{equation}
with the {\it similarity variable}
\begin{equation}
r={x/ t^{\nu  }} >0\,,
\label{eqno(14)}
\end{equation}
with $x > 0, \,t > 0, \, r>0.$
Above  the {\it auxiliary functions} $F_\nu(r)$ and $M_\nu(r)$  are Wright functions (of the second type)  defined 
in the whole complex domain  $z\in \CC$  for $ 0<\nu  <1$ as follows:
\begin{equation}
 F_\nu(z)
  = 
 {\ds  {1\over 2\pi i}\,\int_{Ha}	\!\!  \e^{\ds \, \sigma -z\sigma ^\nu } \,d\sigma}
  =
 {\ds \sum_{n=1}^{\infty}{(-z)^n\over n!\, \Gamma(-\nu  n)}}\,,
       \label{eqno(16)}
   \end{equation}
\begin{equation}
 M_\nu(z) 
 =  
  {\ds {1\over 2\pi i}\,\int_{Ha}   \!\! \e^{\ds \,\sigma -z\sigma ^\nu } \,
   {d\sigma\over \sigma ^{1-\nu  }}} 
   =
 {\ds \sum_{n=0}^{\infty}
 {(-z)^n \over n!\, \Gamma[-\nu  n + (1-\nu  )]}}\,,   
  \label{eqno(17)}
  \end{equation}
where {\it Ha} denotes the Hankel path properly defined for the representation of  the reciprocal  of the  Gamma function.
%%%%%%%%%%\vfill\eject%%%%%%%%%%%%%
%\subsection{ Plots of the Green functions}
%% The Evolution of a Box Function in the Cauchy Problem}

\begin{figure}
\begin{center}
 \includegraphics[height=5cm]{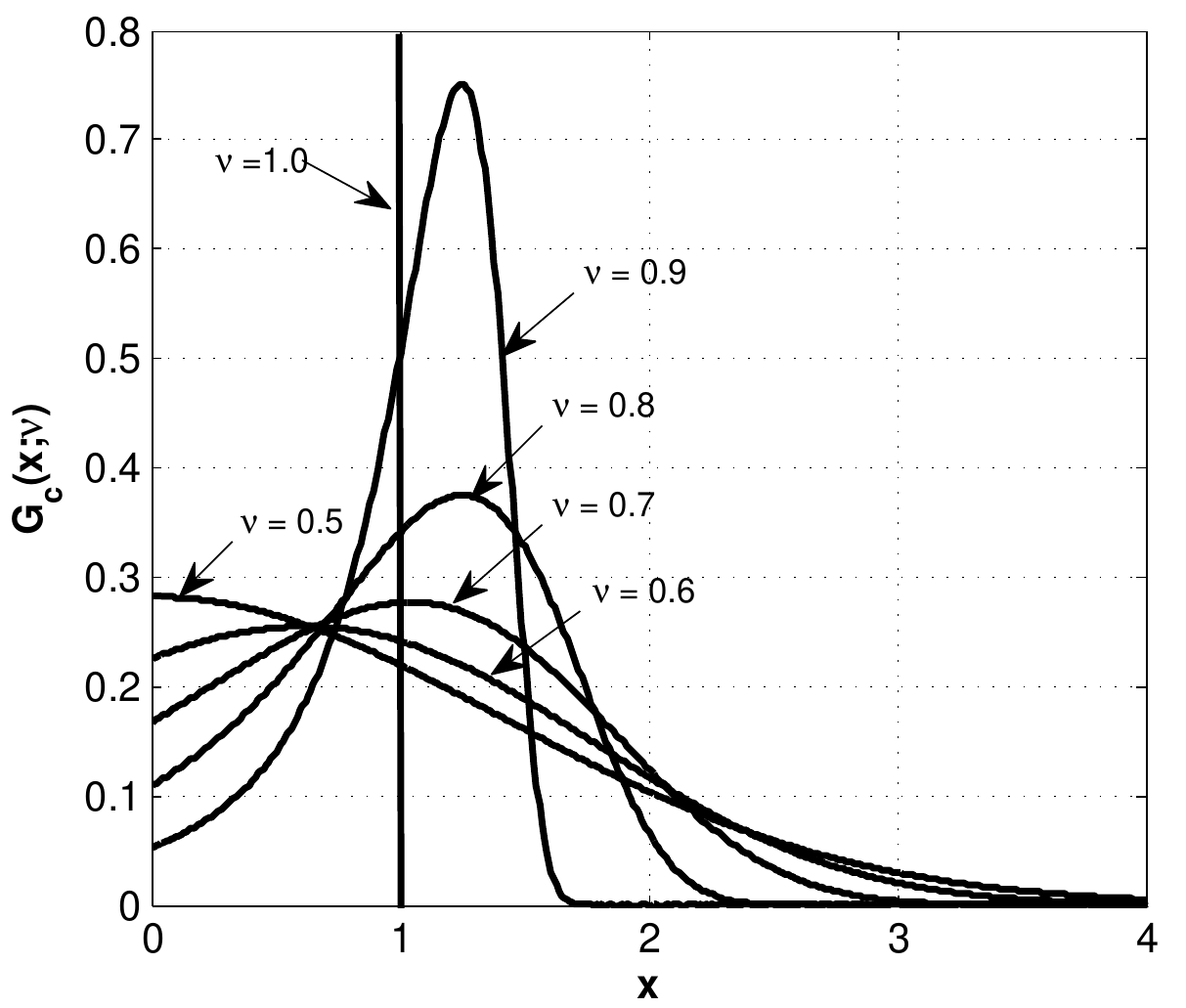} 
\includegraphics[height=5cm]{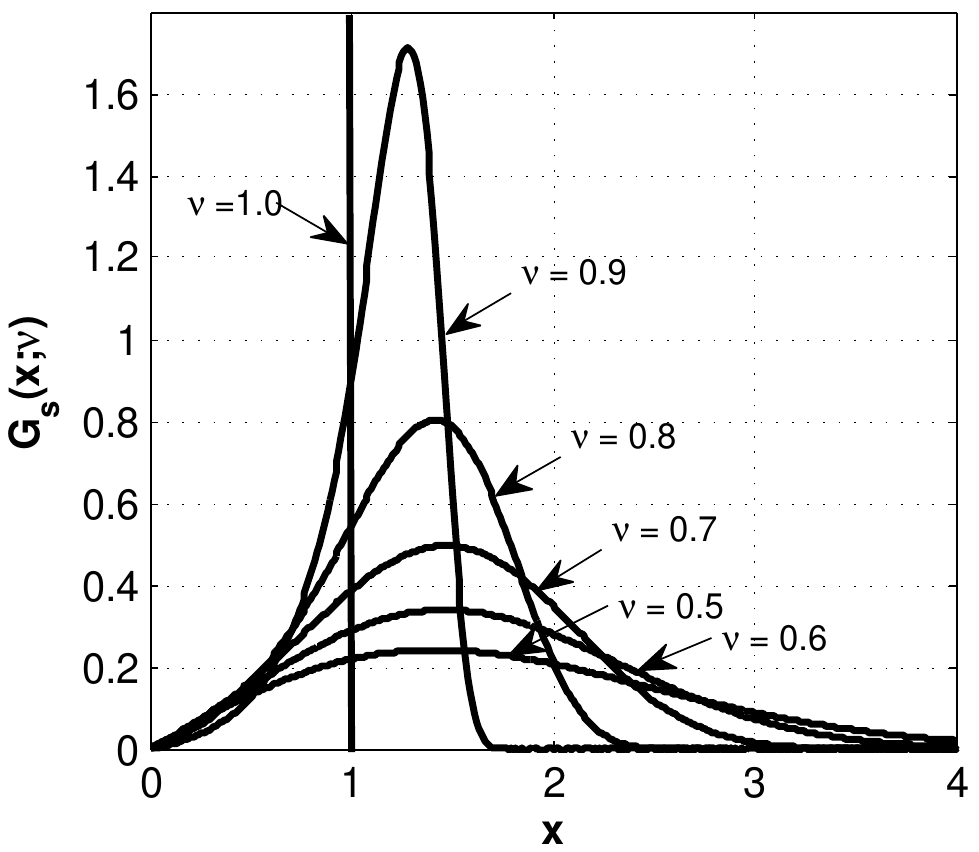}    
\caption{The Green function $\Gc(x;\nu):=\Gc(x,1;\nu)$ (left) and the Green function $\Gs(x;\nu):=\Gs(x,1;\nu)$ (right). Plots for several different values of $\nu$.}
\label{fig:Green}
\end{center}
\end{figure}

In Fig. \ref{fig:Green},  some plots of the Green functions
for both the Cauchy and signaling problems
for different values of $\nu$ and for the fixed time instant $t=1$ are presented. For numerical methods that were used to produce these plots we refer the interested reader to \cite{Gor02}, \cite{Luchko_2008}, \cite{Luchko-Mainardi_2013}, and \cite{LMP_2013}.  
%%%%%%
%%%%%%%%%%%%%%%%%%%%%%%5
%\vskip 0.5 truecm

For more details regarding the auxiliary functions and their properties we refer the interested reader to the appendix F of the book by Mainardi \cite{Mainardi_BOOK2010}.
%We also recall that the auxiliary function $M$ is known as the {\it M-Wright function} and sometimes referred to as Mainardi function, see e.g the book by 
%Podlubny \cite{Podlubny_BOOK1999}.
%The function $W_{\lambda ,\mu }$ was named after  Wright in honor  of the British mathematician E.M. Wright who formerly introduced this function in 1933 \cite{Wright_1933} (first kind) and in 1940\cite{Wright_1940}  (second kind). Unfortunately, in the Bateman project handbook \cite{BATEMAN_1955}, the authors considered only the Wright functions of the first kind. Some noteworthy properties of the Wright functions of the second kind are outlined in the 1970 paper by Stankovi{\'{c}} \cite{Stankovic_1970}
% and in the papers  by Mainardi and Tomirotti \cite{Mainardi-Tomirotti_TMSF94}
%and by Gorenflo, Luchko and Mainardi \cite{GoLuMa_1999,GoLuMa_2000}.

\section{Maximum locations of the Green functions}
In this section, we deal with the maximum locations, maximum values, and propagation velocities of the maximum points of the Green functions $\Gc$ and $\Gs$. Our analysis follows the results presented in Luchko, Mainardi and Povstenko \cite{LMP_2013} for the Cauchy problem and in Luchko, Mainardi \cite{Luchko-Mainardi_2013} for the signaling problem and is restricted to the most interesting and important case of $1/2 < \nu <1$. 

\subsection{The Cauchy problem}
In the paper by Fujita \cite{Fujita_1990a}, a probabilistic proof of the fact that the Green function 
$\Gc(x,t;\nu)$ of the Cauchy problem attains its maximum at the points 
\begin{equation}
\label{max}
x_{*}(t) = \pm c_\nu t^{\nu},\ \nu = \beta/2
\end{equation}
for each $t>0$, where $c_\nu >0$ is a constant determined 
by $\nu,\ 1/2<\nu
 <1$, has been presented for the first time.  In Luchko, Mainardi and Povstenko \cite{LMP_2013}, an analytical proof of this relation was given. 

As mentioned in Fujita \cite{Fujita_1990a}, the maximum point of the Green function 
$\Gc(x,t;\nu)$ propagates for $t>0$ with a finite velocity $\V_c(t,\nu)$ that is determined by
\begin{equation}
\label{speed1}
\V_c(t,\nu):=x_{*}^{\prime}(t) = \nu c_\nu t^{\nu-1}.
\end{equation}
This formula shows that for every $\nu,\ 1/2 < \nu < 1$ 
the propagation velocity of the maximum point of the Green function $\Gc$ 
is a decreasing function in $t$ that varies from $+\infty$ at time $t=0+$ to zero as $t\to +\infty$. 

\begin{figure}
\begin{center}
 \includegraphics[height=5cm]{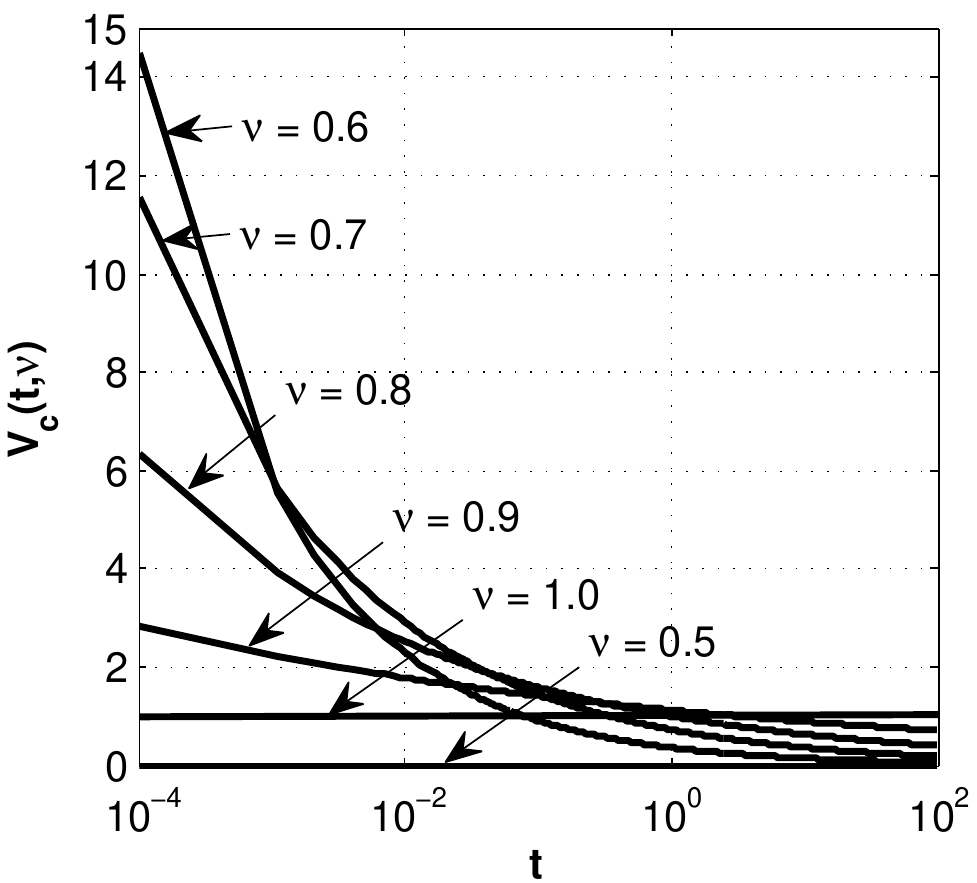} 
\includegraphics[height=5cm]{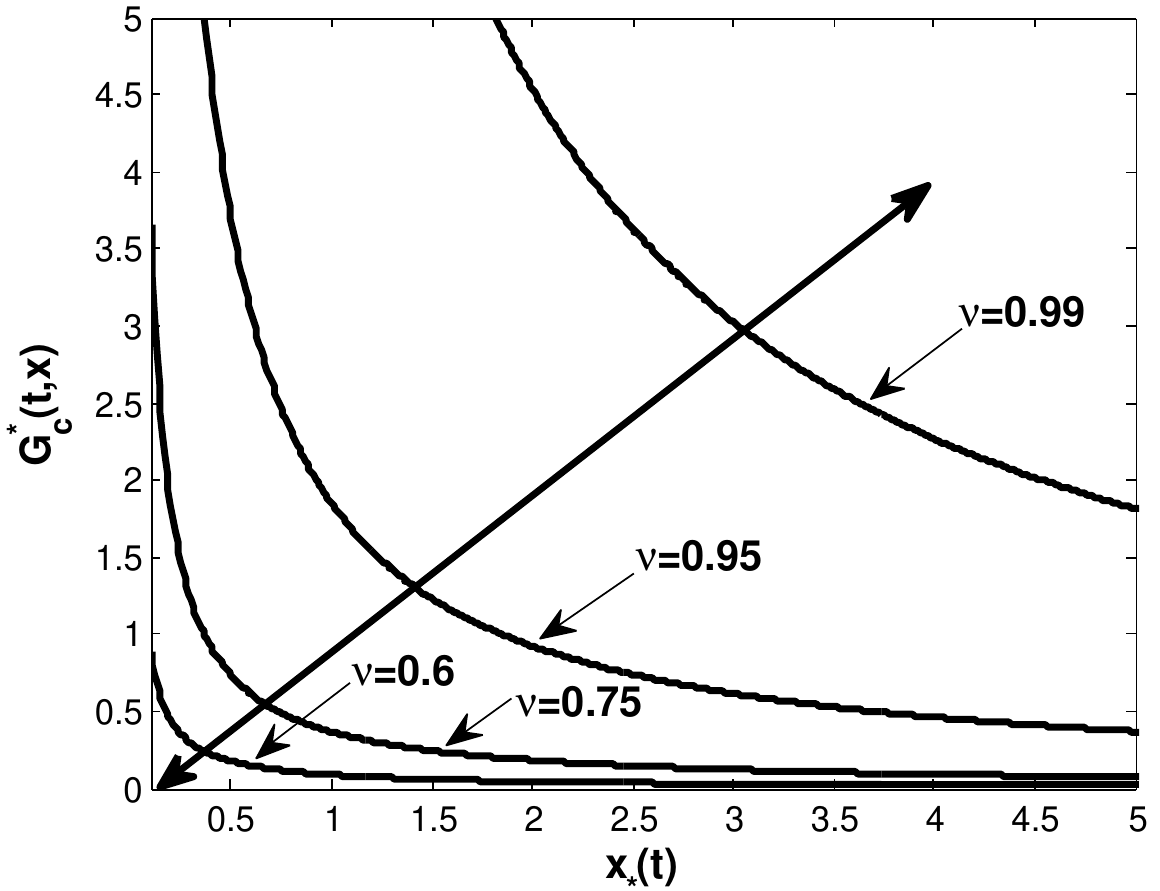}
\caption{Propagation velocity $\V_c(t,\nu)$ of the maximum point of $\Gc$ in the log-lin scale (left) and the maximum locations and maximum values of $\Gc(x,t;\nu)$ (right)}  
%Plots of the parametric curve ($x_{*}(t),\ \Gc^{*}(t;\nu)$),\ $0<t<\infty$ 
%for different values of $\nu$}
\label{fig:speed}
\end{center}
\end{figure}

For $\nu = 1/2$ (diffusion), the propagation velocity is equal to zero because of  $c_{1/2} = 0$ 
whereas for $\nu = 1$ (wave propagation) it remains constant and is equal to $c_{1} = 1$.  
In Fig. \ref{fig:speed}, some plots of the propagation velocity of the maximum point of the Green function $\Gc$ 
are given for different values of $\nu$. 

The maximum value $\Gc^{*}(t;\nu)$ of $\Gc(x,t;\nu)$   is given by the formula
\begin{equation}
\label{MaxVal2}
\Gc^{*}(t;\nu)= m_\nu t^{-\nu},\ \ \ 
%\end{equation}
%\begin{equation}
%\label{MaxVal3}
m_\nu= \frac{1}{2}M_\nu(c_\nu) = \frac{1}{\pi}\int_{0}^{\infty}
             E_{2\nu}\left( -\tau^{2} \right)\, \cos (c_\nu \tau)\, d\tau,
\end{equation}
$E_\alpha$ being the Mittag-Leffler function (see \cite{LMP_2013}).

It follows from the relations (\ref{max}) and (\ref{MaxVal2}) (and of course directly from the formula (\ref{eqno(15)}))  that the product 
\begin{equation}
\label{MaxVal4}
\Gc^{*}(t;\nu)\cdot x_{*}(t) = c_\nu\, m_\nu,\q 0<t<\infty,
\end{equation}
is a constant that depends only on $\nu$, 
i.e., that the maximum locations and the corresponding maximum values specify a 
certain hyperbola for a fixed value of $\nu$ and for $0<t<\infty$. 
%This fact easily follows from the scaling property of the Green function $\Gc$. 
%Let us note that the product $\Gc^{*}(t;\nu)\cdot x_{*}(t)$ 
%is equal to zero in the case $\nu = 1/2$ (diffusion equation) 
%because the maximum point is always located at the point $x_{*}=0$ 
%and to infinity in the case $\nu = 1$ (wave equation) 
%because the maximum value is always equal to infinity. 

In Fig. \ref{fig:speed}, we give some plots of the parametric curve 
($x_{*}(t),\ \Gc^{*}(t;\nu)$),\ $0<t<\infty$ that is in fact a hyperbola  for different values of $\nu$.
 The vertex of the hyperbola tends to the point $(0,0)$ as $\nu$ tends to $1/2$ (diffusion equation) 
 and to infinity as $\nu \to 1$ (wave equation).

\begin{figure}
\begin{center}
\includegraphics[height=5cm]{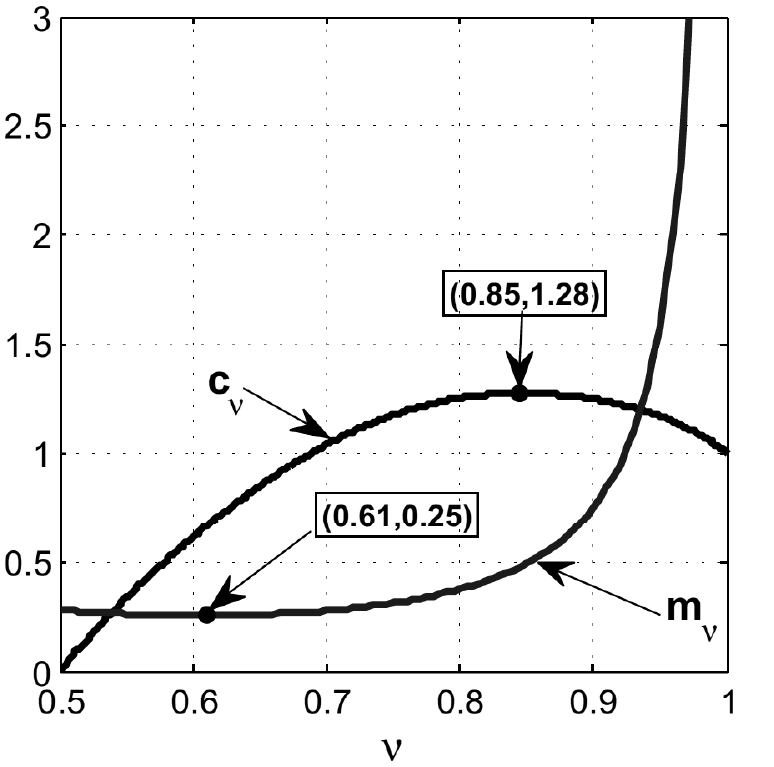} 
\includegraphics[height=5cm]{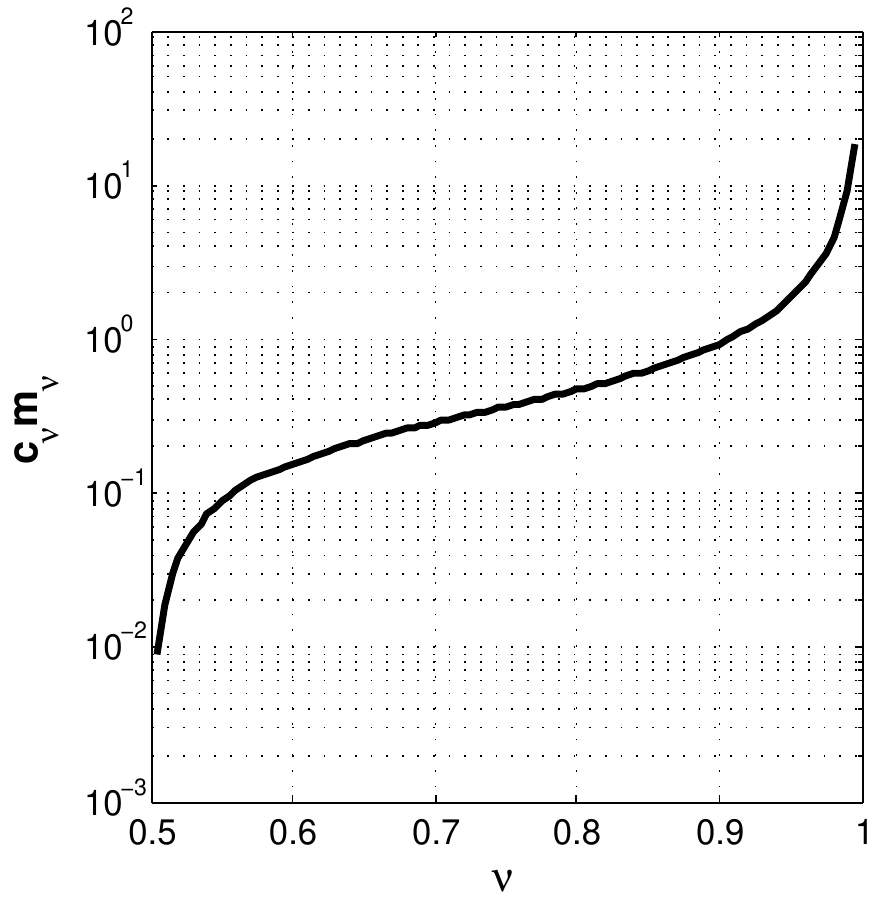}
\caption{Maximum locations and maximum values of the Green function $\Gc(x,1;\nu)$ (left) and their product (right)}
\label{fig:Prod}
\end{center}
\end{figure}

Finally, in Fig. \ref{fig:Prod},  
 the maximum locations $c_\nu$ and the maximum values $m_\nu$ as well as their product  are plotted for $1/2 < \nu <1$.  
%As we have seen above, the constants $c_\nu$ (maximum location of $\Gc(x,1;\nu)$) 
%and $m_\nu$ (maximum value of $\Gc(x,1;\nu)$) 
%are decisive for the behavior of the Green function $\Gc(x,t;\nu)$ for all $t>0$ 
%because the maximum locations and values of this function can be determined via these 
%constants for any time point $t>0$ (see the formulas (\ref{max}) and (\ref{MaxVal2})). 
As expected, the product $c_\nu\, m_\nu$ is a  monotonically increasing  function 
that takes values between $0$ (diffusion equation) and  $+\infty$ (wave equation). 

 \subsection{The Signaling problem}
In this subsection, we present some results regarding the maximum location and the maximum value of the Green function $\Gs$ for the signaling problem as a function in the spatial variable $x$ for the fixed values of $\nu$ and $t$ (see \cite{Luchko-Mainardi_2013} for more details). 

\begin{figure}
\begin{center}
\includegraphics[height=5cm]{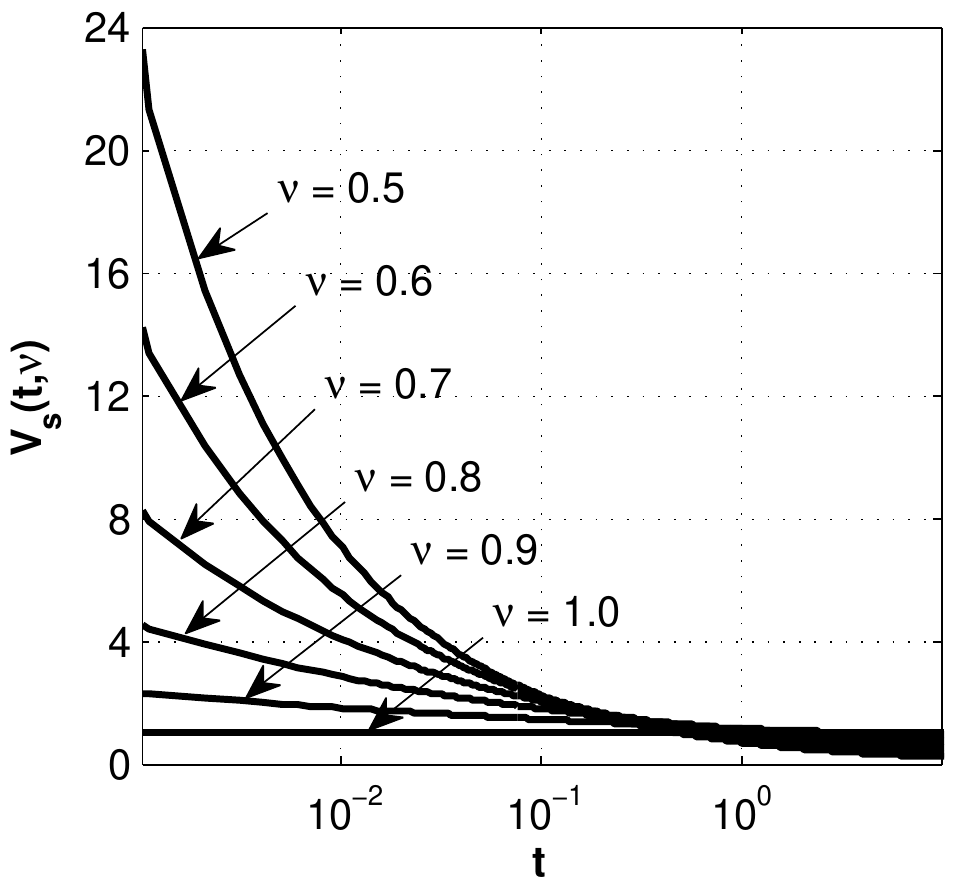} 
\includegraphics[height=5cm]{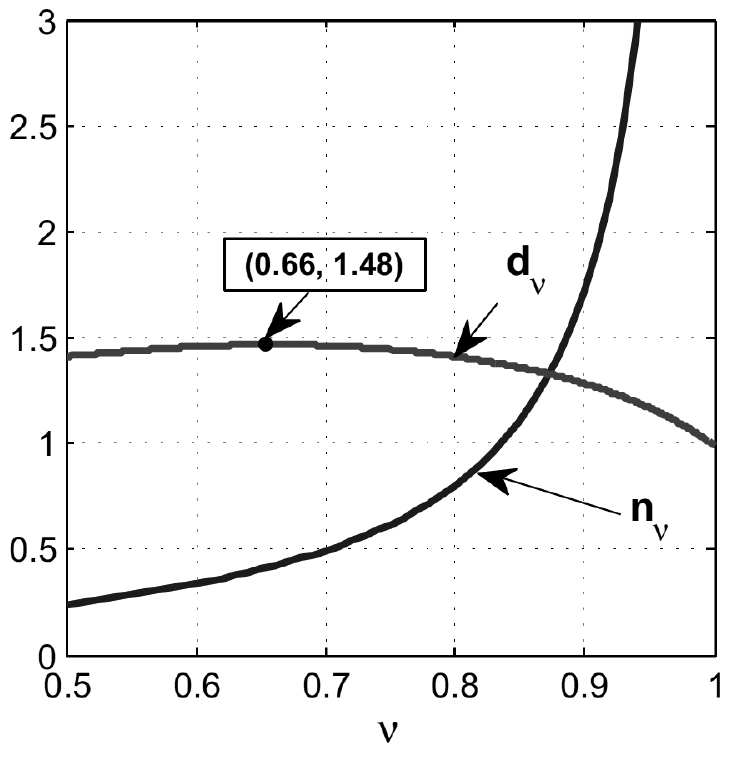}
\caption{Propagation velocity $\V_s(t,\nu)$ of the maximum point of $\Gs$ (left) and the maximum locations and maximum values of the Green function $\Gs(x,1;\nu)$ for $1/2 \le \nu \le 1$ (right)}
\label{fig:speed_s}
\end{center}
\end{figure}

As in the previous subsection,  the formula (\ref{eqno(15)}) is employed to get the following representation for the maximum location $x_{*}=x_{*}(t,\nu)$ of the Green function $\Gs(x,t;\nu)$: 
\begin{equation}
\label{tmax}
x_{*}(t,\nu) = D\, t^{{\nu}}\ \  \mbox{with}\ \  
D= x_{*}(1,\nu)=d_\nu.
\end{equation}
Having determined the  maximum location  of $\Gs$, we can now calculate the propagation velocity $\V_s(t,\nu)$ of the maximum point. If follows from (\ref{tmax}) that
\begin{equation}
\label{vel_max}
\V_s(t,\nu)=\frac{dx_{*}}{dt} ={\nu}{d_\nu}\, t^{{\nu}-1}. 
\end{equation}
As we see, the propagation velocity $\V_s(t,\nu)$ of the maximum point of the Green function $\Gs$ is described by a formula of the same type as the one for the Green function $\Gc$ (see the formula (\ref{speed1})). Its qualitative behavior is therefore very similar to that presented in 
Fig. \ref{fig:speed} as one can see on the plots of Fig.  \ref{fig:speed_s}. The only essential difference between the plots of Figs. \ref{fig:speed} and \ref{fig:speed_s} is in the curve for $\nu = 1/2$ (diffusion equation). Whereas the maximum location of the Green function for the Cauchy problem for the diffusion equation  does not move with the time ($c_{1/2} = 0$), its velocity for the signaling problem is  equal to $\frac{\sqrt{2}}{2}\, \frac{1}{\sqrt{t}}$ (see the formula (\ref{vel_max})).   

The maximum value $\Gs^{*}(t;\nu)$ of $\Gs(x,t;\nu)$ in dependence of $t$ and $\nu$ is given by the following formula: 
\begin{equation}
\label{MaxVal_s}            
 \Gs^{*}(t;\nu) = \frac{n_\nu}{t}\ \ \mbox{with} \ \ 
 n_\nu=    F_\nu(d_\nu) =        
     \frac{2}{\pi}\int_{0}^{\infty}\tau
             E_{2\nu,2\nu}\left( -\tau^{2} \right)\, \sin (d_\nu \tau)\, d\tau,
\end{equation}
$E_{\alpha,\beta}$ being the generalized Mittag-Leffler function (see  \cite{Luchko-Mainardi_2013}).
%Of course, we can also employ the representation (\ref{eqno(15)}) of the Green function $\Gs$ to get the formula  
%\begin{equation}
%\label{recip-1}
% \Gs^{*}(t;\nu) = \Gs(x_{*}(t),t;\nu  ) 
%    = \frac{1}{t} F_\nu(x_{*}(t)/t^{\nu })= \frac{1}{t} F_\nu(d_\nu).
%\end{equation}
%Comparing (\ref{MaxVal_s}) and  (\ref{recip-1}), we get the relation
%$$
%n_\nu = F_\nu(d_\nu)
%$$
%that shows that $d_\nu$ is the maximum location of the function $F_\nu$ and $n_\nu$ is its maximum value. 
%
%Evidently, $n_\nu$ depends just on $\nu$ but not on the time variable $t$. 
%For $\nu = 1/2$, the maximum value $\Gs^{*}(t;\nu)$ can be easily determined from the representation (\ref{tmax-1/2}):
%\begin{equation}
%\label{Vmax-1/2}
%\Gs^{*}(t;1/2) = \frac{1}{\sqrt{2\pi\e}} t^{-1},
%\end{equation}
%i.e., we get
%\begin{equation}
%\label{n-1/2}
%n_{\frac{1}{2}} = \frac{1}{\sqrt{2\pi\e}} \approx  0.2420.
%\end{equation}
%If $\nu = 1$, the maximum value of  $\Gs$  is evidently equal to $+\infty$, i.e.,
%\begin{equation}
%\label{n-1}
%n_{1}= +\infty.
%\end{equation}
The maximum locations and the maximum values of the Green function $\Gs(x,1;\nu)$ for the intermediate values of $\nu,\ 1/2 < \nu < 1$ are presented in Fig. \ref{fig:speed_s}. 
%It is interesting to note that  the function $d_\nu$ is not monotone on the interval $\nu \in [1/2, 1]$. 
%As we can see on the plot, the function $d_\nu = d_\nu(\nu),\ 1/2 \le \nu \le 1$ has a maximum located at the point $\nu \approx 0.66$. The maximum value of the maximum location is approximately equal to $1.48$.  As to the maximum value of the Green function $\Gs(x,t;\nu)$ at the time instant $t=1$, it is a monotonically increasing function in $\nu$ that starts with $\frac{1}{\sqrt{2\pi\e}} \approx  0.24$ at the point $\nu = 1/2$ and tends to $+\infty$ as $\nu \to 1$. 

Let us finally note that it follows from the formulas (\ref{tmax}) and (\ref{MaxVal_s}) (and of course directly from the formula (\ref{eqno(15)})) that the product $p_\nu$ of the maximum location and the maximum value of the Green function $\Gs(x,t;\nu)$ is time-dependent
\begin{equation}
\label{pnu}
p_\nu = p_nu(t)= \Gs^{*}(t;\nu)\cdot x_{*}(t) = d_\nu t^{{\nu}}\, n_\nu t^{-1} = d_\nu n_\nu t^{{\nu-1}} 
\end{equation}
and follows the formula of the same type as the one for the propagation velocities of the maximum locations for the Green functions $\Gs$ and $\Gc$.

\section {Centers of gravity and medians of \\ the Green functions $\Gc$ and $\Gs$}
In this section, some new results regarding locations and velocities of the centers of gravity and medians of the Green functions both for the Cauchy and the signaling problems are presented. The key role in all calculations is played by the formula (\ref{eqno(15)}) that connects the Green functions $\Gc$ and $\Gs$ with the special functions of the Wright type. For the readers convenience we list here some formulas for the Mainardi function $M_\nu$ that are used in the further discussions. 

The asymptotics of $M_\nu$ is described by the following formula (see e.g. \cite{Mainardi-Luchko-Pagnini_FCAA2001}): 
\begin{equation}
\label{asymp}
\begin{array} {ll} 
M_\nu (r)  \sim  A_0  \,
     Y^{\, \nu - 1/2 } \, \exp \,( - Y)  \,,
\q  r \to \infty\,,
\\ 
{\ds A_0 = {1\over   \sqrt{2\pi}\,(1-\nu )^\nu \, \nu ^{2\nu-1}} \,,\ q 
Y = (1-\nu )\, (\nu^\nu\, r)^{1/(1-\nu )}.} 
\end{array}
\end{equation}
As $r\to 0$, $M_\nu (r)$ evidently tends to $1/\Gamma(1-\nu)$. 

The known Mellin transform of the Wright function 
\begin{equation}
\label{(1.67-wright)}
\begin{array}{ll}
{\mathcal M}\{  W_{\lambda,\mu}(-x);s\}  =
\frac{\Gamma(s)}
{\Gamma(\mu -\lambda s)}
\\
0 < \Re(s),\ \lambda < 1  \; \mbox{or}  \ ;
 0 < \Re(s) < \Re(\mu)/2 -1/4, \ \lambda=1
 \end{array}
 \end{equation}
leads to the following formula for the Mellin transform of the Mainardi function $M_\nu$:
\begin{equation}
\label{Mellin_M}
{\mathcal M}\{  M_{\nu}(x);s\} =
\frac{\Gamma(s)}
{\Gamma(1-\nu +\nu s)},\ \ 
0 < \Re(s),\ -1 < \nu.
\end{equation}
Let us note that the {\it  Mellin integral transform} of a sufficiently well-behaved
function $f$ is defined as 
\begin{equation}
\label{(1.49)}
{\mathcal M} \{ f(x);s\}=f^*(s)=\int_0^{+\infty}f(x)x^{s-1}dx.
\end{equation}

\subsection{Center of gravity of the Green function  for the Cauchy problem}
Let us start with calculation of the center of gravity of the Green function $\Gc$. Because  $\Gc$ is an even function, we restrict our attention to the function $\Gc(r,t;\nu) = \Gc(|x|,t;\nu),\ r = |x|\ge 0$. The location $r_\nu^{c}(t)$ of the center of gravity of $\Gc(r,t;\nu)$  is defined  by the formula  
\begin{equation}
\label{grav}
r_\nu^{c}(t) = \frac{\int_0^\infty r\, \Gc(r,t;\nu)\, dr}
{\int_0^\infty \Gc(r,t;\nu)\, dr}.
\end{equation}
Using the formulas (\ref{eqno(15)}) and (\ref{Mellin_M}) and a linear variables substitution, we can calculate both integrals in (\ref{grav}) in explicit form. For the first integral we get
$$
\begin{array}{ll}
{\ds \int_0^\infty r\, \Gc(r,t;\nu)\, dr} 
&= {\ds \int_0^\infty r\, \frac{1}{2t^\nu}M_\nu (r/t^\nu)\, dr = \frac{t^\nu}{2} \int_0^\infty u\, M_\nu (u)\, du} 
\\ \\
&= {\ds \frac{t^\nu}{2} \frac{\Gamma(2)}{\Gamma(1-\nu +2\nu)} =
\frac{t^\nu}{2 \Gamma(1+\nu)}\,.}
\end{array}
$$
The second integral in (\ref{grav}) is evidently equal to $1/2$ because $\Gc(x,t;\nu)$ is a probability density function (pdf)  in $x$ evolving in time that is an even function (see e.g. \cite{Mainardi-Luchko-Pagnini_FCAA2001}). Of course, this integral can be calculated explicitly following the same method we applied for the first integral. 

The final formula for the location of the center of gravity  of the Green function $\Gc$ of the Cauchy  problem is as follows:
\begin{equation}
\label{grav-final}
r_\nu^{c}(t) = g_c(\nu)\,t^\nu, \ \ g_c(\nu)=\frac{1}{\Gamma(1+\nu)}.
\end{equation}
For the diffusion equation with $\nu = 1/2$, the location $r_{1/2}^{c}(t)$  of the center of gravity for a fixed $t$ is at the point $\frac{2\sqrt{t}}{\sqrt{\pi}}$ whereas for the wave equation we get as expected $r_{1}^{c}(t) = t$. 

\begin{figure}
\begin{center}
\includegraphics[height=5cm]{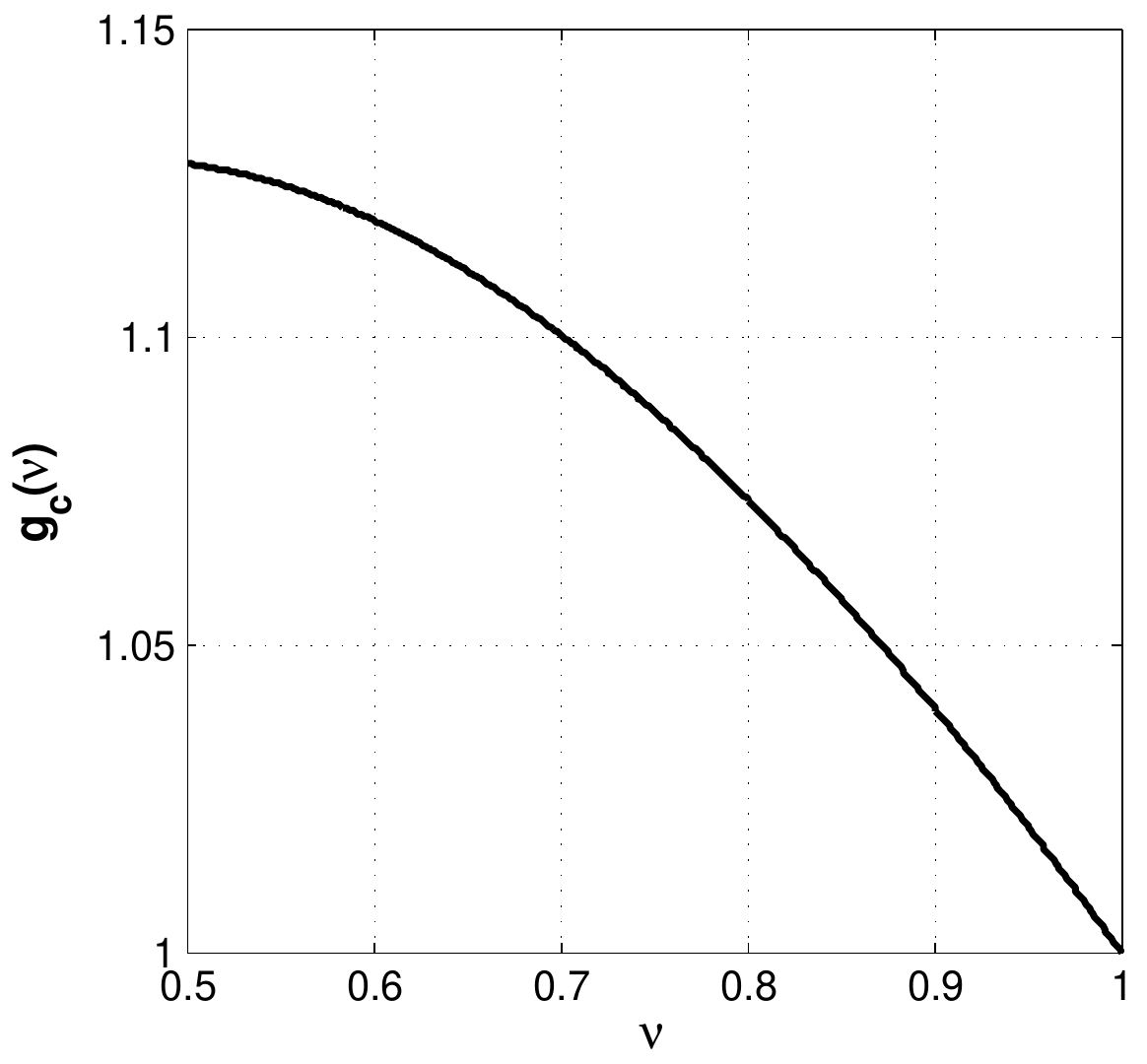} 
\includegraphics[height=5cm]{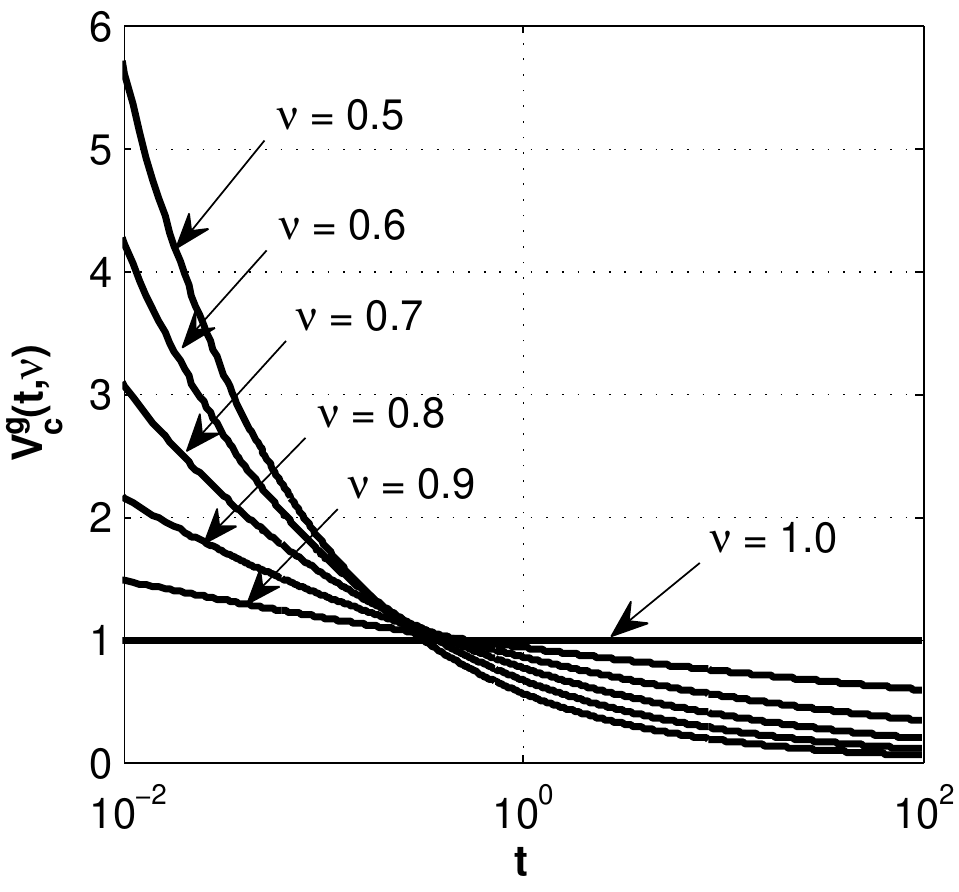}
\caption{Location of the center of gravity of the Green function $\Gc(|x|,1;\nu)$ (left) and propagation velocity $\V_c^{g}(t,\nu)$ of the center of gravity of $\Gc$ (right)}
\label{fig:gc}
\end{center}
\end{figure}

As we can see from the formula (\ref{grav-final}), the location of the center of gravity is a power function in $t$ with the coefficient $g_c(\nu)$ that depends on $\nu$. On the other hand, $g_c(\nu)$ describes the location of the center of gravity of $\Gc$ at the time instant $t = 1$. A plot of the function $g_c(\nu),\ 1/2\le \nu \le 1$ is presented in Fig. \ref{fig:gc}. Because $\Gamma(1+\nu)$ is a monotonically increasing function for $1/2 \le \nu \le 1$, the function $g_c(\nu)$ monotonically decreases from $\frac{2}{\sqrt{\pi}} \approx  1.1284$  (diffusion equation) to $1$ (wave equation). 

The velocity $\V_c^{g}(t,\nu)$ of the center of gravity of $\Gc$ is given by the formula
\begin{equation}
\label{vel_g_c}
\V_c^{g}(t,\nu) = \frac{d}{dt} r_\nu^{c}(t) = \frac{t^{\nu-1}}{\Gamma(\nu)},
\end{equation}
so that again we obtain a formula of the same type as the one for velocities of the maximum locations of the Green functions (with a different coefficient). The plots of the velocity $\V_c^{g}(t,\nu)$ for different values of $\nu$ look like the ones presented in Fig. \ref{fig:speed_s} and are given in Fig. \ref{fig:gc}. 

\subsection{Center of gravity of the Green function for the Signaling problem} 
We consider now the location of the center of gravity of the Green function $\Gs$ of the signaling problem.
It is known (see e.g. \cite{Luchko-Mainardi_2013}) that for a fixed $x>0$ and 
for a fixed $\nu,\ 1/2\le \nu <1$, the Green function $\Gs(x,t;\nu)$ is a one-sided {\it stable}
probability density function (pdf) of the time variable $t>0$. 
A prominent example is the function $\Gs(x,t;1/2)$ that for $x=1$ is called the L\'evy-Smirnov pdf. It follows from the formulas (\ref{eqno(15)}) and (\ref{asymp}) that the location of the gravity center of $\Gs$ with respect to the time variable $t$ 
is in infinity for all $x>0$. That is why we consider the location of the gravity center of $\Gs$ with respect to the spatial variable $x$ that is defined by the formula
\begin{equation}
\label{grav_s}
r_\nu^{s}(t) = \frac{\int_0^\infty r\, \Gs(r,t;\nu)\, dr}
{\int_0^\infty \Gs(r,t;\nu)\, dr}.
\end{equation}
To evaluate the integrals in the formula (\ref{grav_s}), we again employ the relation (\ref{eqno(15)}), a linear variables substitution in the integrals, and the formula (\ref{Mellin_M}) for the Mellin transform of the Mainardi function $M_\nu$ and thus get the following results:
$$
\int_0^\infty \Gs(r,t;\nu)\, dr = \nu t^{\nu -1} \int_0^\infty u\, M_\nu(u)\, du = \frac{\nu }{\Gamma(1+\nu)} t^{\nu -1},
$$
$$
\int_0^\infty r\, \Gs(r,t;\nu)\, dr =  \nu t^{2\nu -1} \int_0^\infty u^2\, M_\nu(u)\, du = \frac{2\nu }{\Gamma(1+2\nu)} t^{2\nu -1}.
$$
Using the duplication and the reduction formulas for the Gamma function, we obtain the following formula for the location of the center of gravity of the Green function $\Gs$: 
\begin{equation}
\label{grav_s_f}
r_\nu^{s}(t)
 = \frac{\frac{2\nu }{\Gamma(1+2\nu)} t^{2\nu -1}}
{\frac{\nu }{\Gamma(1+\nu)} t^{\nu -1}} 
= \frac{\Gamma(\nu)}{\Gamma(2\nu)} \, t^\nu 
=  \frac{\sqrt{\pi}\, 2^{1-2\nu}}{\Gamma\left(\nu + \frac{1}{2}\right)}\, t^\nu.
\end{equation}
The two known particular cases of this formula are $r_{1/2}^{s}(t) =\sqrt{\pi}\, t^{1/2}$ (diffusion equation) and $r_{1}^{s}(t) =t$ (wave equation). 

\begin{figure}
\begin{center}
\includegraphics[height=5cm]{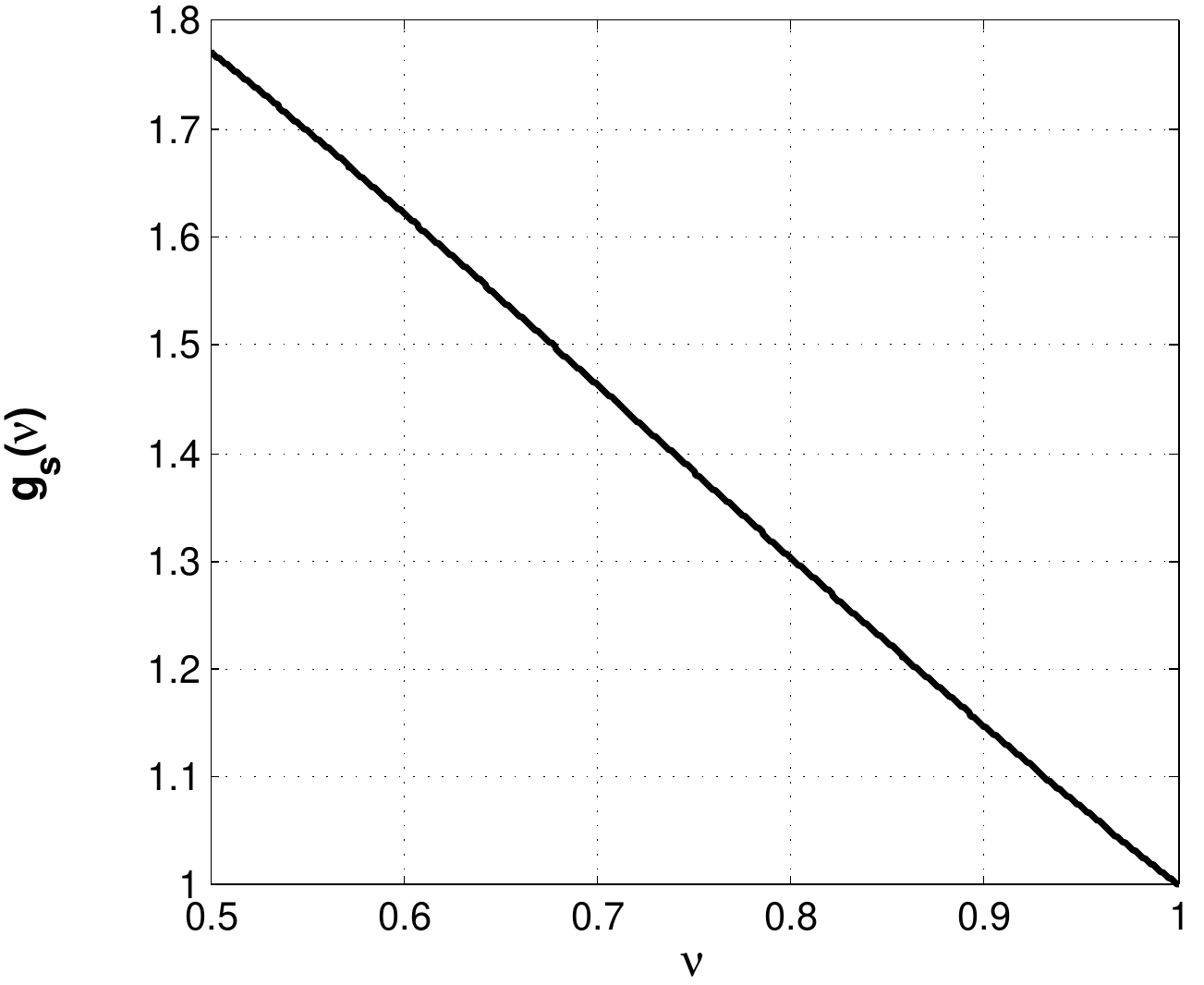} 
\includegraphics[height=5cm]{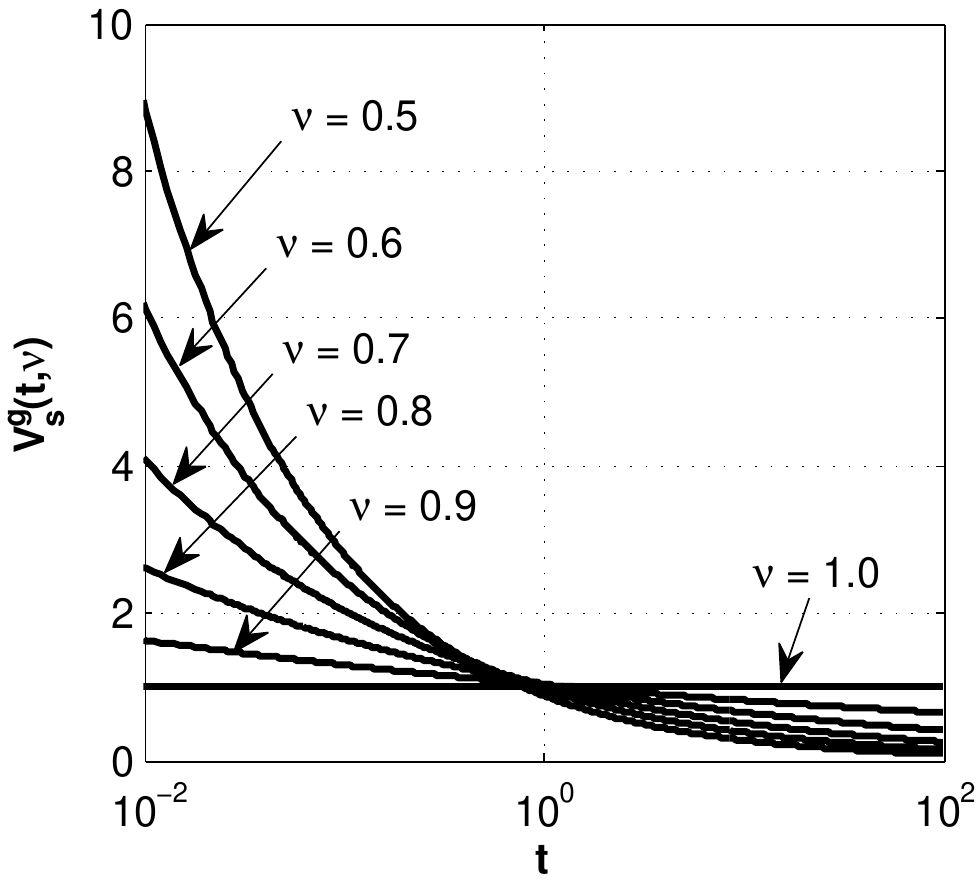}
\caption{Location of the center of gravity of the Green function $\Gs(x,1;\nu)$ (left) and propagation velocity $\V_s^{g}(t,\nu)$ of the center of gravity of $\Gs$ (right)}
\label{fig:gs}
\end{center}
\end{figure}

At the time instant $t=1$, we get the relation
\begin{equation}
\label{gs}
g_s(\nu):= r_\nu^{s}(1)  = \frac{\sqrt{\pi}\, 2^{1-2\nu}}{\Gamma\left(\nu + \frac{1}{2}\right)}.
\end{equation}
A plot of the function $g_s(\nu)$ is presented in the Fig \ref{fig:gs}. As we can see, the function $g_s(\nu)$ monotonically decreases from the value $\sqrt{\pi}\approx 1.7725 $ at the point $\nu = 1/2$ (diffusion equation) to the value 1 at the point $\nu = 1$ (wave equation). Surprisingly, the plot of $g_s(\nu)$ is very similar to a straight line, but of course $g_s(\nu)$ is not a linear function on the interval $1/2\le \nu \le 1$.

The velocity $\V_s^{g}(t,\nu)$ of the center of gravity of $\Gs$ is calculated via the formula
\begin{equation}
\label{vel_g_s}
\V_s^{g}(t,\nu) = \frac{d}{dt} r_\nu^{s}(t) = \frac{\sqrt{\pi}\, 2^{1-2\nu}\,  \nu}{\Gamma\left(\nu + \frac{1}{2}\right)}\, t^{\nu-1}.
\end{equation}
Some plots of the velocity $\V_s^{g}(t,\nu)$ for different values of $\nu$ are presented in Fig. \ref{fig:gs} (of course, they are very similar to those shown in  Fig. \ref{fig:gc}).

\subsection{Medians of the Green functions for the Cauchy and Signaling problems}
Finally, we consider the locations of the medians of the Green functions $\Gc$ and $\Gs$. 

The location $x = x_m^c$ of the median of the Green function $\Gc(r,t;\nu),\ r\ge 0$ can be determined from the equation
$$
\int_0^{x_m^c} \Gc(r,t;\nu)\, dr = \frac{1}{2} \int_0^{\infty} \Gc(r,t;\nu)\, dr = \frac{1}{4}.
$$
Because
$$
\int_0^{x_m^c} \Gc(r,t;\nu)\, dr = \int_0^{x_m^c} \frac{1}{2t^\nu}M_\nu(r/t^\nu)\, dr = \int_0^{x_m^c/t^\nu} \frac{1}{2}M_\nu(u)\, du,
$$
the location of the median can be determined from the equation
\begin{equation}
\label{eq_m_c}
\int_0^{x_m^c/t^\nu} \frac{1}{2}M_\nu(u)\, du = \frac{1}{4}.
\end{equation}
The integral $\int_0^{x} M_\nu(u)\, du$ monotonically increases from 0 to 1 as $x$ varies from 0 to $+\infty$ and thus the equation (\ref{eq_m_c}) has a unique solution 
that of course depends on $\nu$ and is denoted by $m_c(\nu)$. 
Then we first get the relation
$$
x_m^c/t^\nu = m_c(\nu)
$$
and then the formula
\begin{equation}
\label{med_g_c}
x_m^c = m_c(\nu)t^\nu 
\end{equation}
for the location of the median. Once again, we see that the location of the median is a power function in $t$ with the exponent $\nu$ and the coefficient $m_c(\nu)$ that corresponds to the location of the median at the time instant $t = 1$. 
%This time, no explicit formula for the coefficient $m_c(\nu)$ was derived, but we can evaluate this coefficient numerically as solution to the equation
%$$ 
%\int_0^{m_c(\nu)} M_\nu(u)\, du = \frac{1}{2}.
%$$
%A plot of the location $m_c(\nu)$  of the median of the Green function $\Gc$ at the time instant $t = 1$ is presented in Figure \ref{med_c}. 

For the fixed $x >0$ and $\nu,\ 1/2 \le \nu \le 1$, the location $t = t_m^s$ of the median of the Green function $\Gs(x,t;\nu)$ of the signaling problem for the fixed $x$ and $\nu$ is determined from the equation
$$
\int_0^{t_m^s} \Gs(x,t;\nu)\, dt =  \frac{1}{2}.
$$
Employing the same arguments as in the case of the Green function $\Gc$, we arrive at the equation
\begin{equation}
\label{eq_m_s}
\int_{x/(t_m^s)^\nu}^{+\infty} M_\nu(u)\, du = \frac{1}{2}
\end{equation}
that evidently has the same solution $m_c(\nu)$ as the equation (\ref{eq_m_c}) because of the fact that the integral $\int_0^{x} M_\nu(u)\, du$ monotonically increases from 0 to 1 as $x$ varies from 0 to $+\infty$. 
It follows from (\ref{eq_m_s}) that
$$
x/(t_m^s)^\nu = m_c(\nu)
$$
and 
\begin{equation}
\label{med_g_s}
t_m^s = m_s(\nu)x^{1/\nu},\ \  m_s(\nu) = \frac{1}{(m_c(\nu))^{1/\nu}}.
\end{equation} 
%A plot of the location $m_s(\nu)$  of the median of the Green function $\Gs$ at the point $x = 1$ is presented in Figure \ref{med_s}. 

\section{Conclusions and discussions}
In this paper, we deal with some important properties of the Green functions
of the Cauchy and signaling problems for the one-dimensional time-fractional diffusion-wave equation with the constant coefficients.
 It is known that these functions are relevant to characterize the evolution of 
 the pulse-like initial data that appears as an intermediate process between diffusion and wave propagation. 
 Except in the limiting case of the wave equation, the pulses propagate with infinite velocities  that is typical for evolution equations of the parabolic type.  
 These processes are common to refer to as the {\it diffusive waves}. 
 In this paper, we show that the maximum locations and the centers of gravity (in space) of the diffusive waves always propagate with a finite velocity that is determined by a power law in time. 
 The exponent of the power law (related to the order of the  fractional derivative in the time-fractional diffusion-wave equation) is the  same for the Cauchy and signaling problems 
as expected from the similarity properties of the evolution equation.  
 Whereas the location of the maxima and their velocities are determined numerically, the corresponding quantities   for the centers of gravity are obtained analytically in a closed form.   
In the absence of a finite wave-front velocity that is typical for the standard waves, the velocities mentioned above can be interpreted as a sort of characteristic signal velocity
of the diffusive waves and thus worthy to be investigated and calculated. 

\section{Acknowledgments}
The authors appreciate constructive remarks and suggestions of the anonymous referees that helped to improve the manuscript. 

\input{BIB-JVA.tex}

\end{document}

%% file: BIB-JVA.tex
%% BIB-JVA.tex